\documentclass[reqno]{amsart}
\usepackage{amssymb,amscd}

\newcommand{\ddq}{\mathop{\rm wt}\nolimits}	
\newcommand{\ebv}{\epsilon}	

\newcommand{\defterm}[1]{{\it #1}}
\newcommand{\comment}[1]{{}}

\newcommand{\sm}{\setminus}
\newcommand{\sd}{\triangle}
\newcommand{\st}{~\mid~}

\newcommand{\x}{\times}
\newcommand{\xvar}[2]{x_{#1}^{(#2)}}	

\newcommand{\dir}{\mathop{\rm dir}\nolimits}

\newcommand{\Tree}{{\mathop{\rm Tree}\nolimits}}
\newcommand{\indeg}{\mathop{\rm indeg}\nolimits}
\newcommand{\outdeg}{\mathop{\rm outdeg}\nolimits}
\newcommand{\dd}{\mathop{\rm dd}\nolimits}
\newcommand{\id}{\mathop{\rm id}\nolimits}

\newcommand{\QQ} {{\bf Q}}
\newcommand{\ZZ} {{\bf Z}}
\newcommand{\condns}[2]{\substack{#1 \\ #2}}

\newcommand{\mm} {{\mathfrak m}}

\newtheorem{theorem}{Theorem}

\newtheorem{proposition}[theorem]{Proposition}

\newtheorem{lemma}[theorem]{Lemma}

\begin{document}
\title{Factorizations of some weighted spanning tree enumerators}

\author{Jeremy L. Martin}
\address{School of Mathematics\\
University of Minnesota\\
Minneapolis, MN 55455}
\email{martin@math.umn.edu}

\author{Victor Reiner}
\address{School of Mathematics\\
University of Minnesota\\
Minneapolis, MN 55455}
\email{reiner@math.umn.edu}


\thanks{First author supported by NSF Postdoctoral Fellowship.
Second author supported by NSF grant DMS-9877047.}

\begin{abstract}
We give factorizations for weighted spanning tree enumerators of Cartesian
products of complete
graphs, keeping track of fine weights related to degree sequences and edge directions.
Our methods combine Kirchhoff's Matrix-Tree Theorem with the technique
of identification of factors.
\end{abstract}

\maketitle

\section{Introduction}

Cayley's celebrated formula $n^{n-2}$ for the number of spanning trees
in the complete graph $K_n$ has many generalizations (see \cite{Moon}).
Among them is the following well-known factorization for the enumerator of
the spanning trees according to their degree sequence, which is a model
for our results.

\vskip .1in
\noindent
{\bf Cayley-Pr\"ufer Theorem.}
{\it
$$
\sum_{T \in \Tree(K_n)} x^{\deg(T)}
= x_1 x_2 \cdots x_n (x_1 + \cdots + x_n)^{n-2}
$$
where $\Tree(K_n)$ is the set of all spanning trees and
$x^{\deg(T)} := \prod_{i=1}^n x_i^{\deg_T(i)}$.
}
\vskip .1in

Although this is most often deduced from the bijective proof of Cayley's
formula that uses
Pr\"ufer coding (see, e.g., \cite[pp. 4-6]{Moon}), we do not know of such
a bijective proof
for most of
our later results. Section~\ref{model-proof-section} gives a quick proof
(modelling
those that will follow) using a standard weighted version of {\it
Kirchhoff's Matrix-Tree
Theorem}, along with the method of {\it identification of factors}.

We generalize the Cayley-Pr\"ufer Theorem to Cartesian
products of complete graphs $G = K_{n_1} \times \cdots \times K_{n_r}$.
The number of
spanning trees for such product graphs can be computed using Laplacian
eigenvalues.
Section~\ref{Cartesian-products-section} generalizes this calculation to
keep track of the
directions of edges in the tree, as we now explain.  Note that vertices in
$G$
are $r$-tuples $(j_1,\ldots,j_r) \in [n_1] \times \cdots \times
[n_r]$, and each edge connects two such $r$-tuples that differ in only one
coordinate.
Say
that such an edge lies {\it in direction} $i$ if its two endpoints differ
in their $i^{th}$
coordinate.  Given a spanning tree $T$ in $G$, define the
direction monomial
$$
q^{\dir(T)}:= \prod_{i=1}^r q_i^{|\{\text{edges in }T\text{ in direction
}i\}|}.
$$

\begin{theorem}
\label{directions-theorem}
$$
\begin{aligned}
\sum_{T \in \Tree(K_{n_1} \times \cdots \times K_{n_r})}
                        q^{\dir(T)}
&= \frac{1}{n_1 \cdots n_r}\prod_{\emptyset \neq A \subset [r]}
    \left( \sum_{i \in A} q_i n_i \right )^{\prod_{i \in A}(n_i-1)} \\
&= \prod_{i=1}^r q_i^{n_i-1} n_i^{n_i-2}
    \prod_{\condns{A \subset [r]}{|A| \geq 2}}
    \left( \sum_{i \in A} q_i n_i \right )^{\prod_{i \in A}(n_i-1)} .
\end{aligned}
$$
\end{theorem}

One might hope to generalize the previous result by keeping track of edge
directions and
vertex degrees simultaneously.  Empirically, however, such generating
functions do not
appear to factor nicely.  Nevertheless, if one ``decouples'' the vertex
degrees in a certain
way that we now explain, nice factorizations occur.  Create a variable
$\xvar{j}{i}$ for each
pair $(i,j)$ in which $i \in \{1,2,\ldots,r\}$ is a direction and $j$ is in
the range
$1,2,\ldots,n_i$. In other words, there are $r$ sets of variables, with
$n_i$ variables
$\xvar{1}{i},\xvar{2}{i},\ldots,\xvar{n_i}{i}$ in the $i^{th}$ set.  Given
a spanning tree $T$
of $G$, define the decoupled degree monomial
\begin{equation}
\label{decoupled-monomial}
x^{\dd(T)} := \prod_{v=(j_1,\ldots,j_r) \:\in\: [n_1] \times \cdots \times
[n_r]}
            \left( \xvar{j_1}{1} \cdots \xvar{j_r}{r} \right)^{\deg_T(v)}.
\end{equation}

In Section~\ref{Cartesian-products-section},
we prove the following generalization of the Cayley-Pr\"ufer Theorem.

\begin{theorem}
\label{divisibility-theorem}
The spanning tree enumerator
$$
f_{n_1,\ldots,n_r}(q,x):=
\sum_{T \in \Tree(K_{n_1} \times \cdots \times K_{n_r})}
                        q^{\dir(T)} x^{\dd(T)}
$$
is divisible by
	$$q_i^{n_i-1},$$
by
	$$\left( \xvar{j}{i} \right)^{n_1 \cdots n_{i-1} n_{i+1} \cdots
n_r}$$
and by
	$$\left( \xvar{1}{i} + \cdots + \xvar{n_i}{i} \right)^{n_i-2}$$
for each $i \in [r]$ and $j \in [n_i]$.
\end{theorem}

\noindent
{\bf Conjecture.}
{\it
The quotient polynomial
$$
\frac{f_{n_1,\ldots,n_r}(q,x)}
{ \displaystyle
\prod_{i=1}^r q_i^{n_i-1} \,\,
\left( \xvar{1}{i} \cdots \xvar{n_i}{i} \right)^{n_1 \cdots n_{i-1} n_{i+1}
\cdots n_r} \:
\left( \xvar{1}{i} + \cdots + \xvar{n_i}{i} \right)^{n_i-2}
}
$$
in $\ZZ[q_i, \xvar{j}{i}]$ has non-negative coefficients.
}

\vskip .1in

Empirically, this quotient polynomial seems not to factor further in
general, although when
one examines the coefficient of particular ``extreme'' monomials in the
$q_i$, the resulting
polynomial in the $\xvar{j}{i}$ factors nicely. Such nice factorizations
seem to fail for the
coefficients of non-extreme monomials in $q_i$ when there are at least two
$n_i \geq 3$.

Section~\ref{cube-section} shows that when all $n_i=2$, the
spanning tree enumerator $f_{n_1,\ldots,n_r}(q,x)$ factors beautifully.
Here we consider
the Cartesian product
	$$Q_n:=\underbrace{K_2 \times \cdots \times K_2}_{n\text{
times}},$$
which is the $1$-skeleton of the $n$-dimensional cube.  For the sake of a
cleaner
statement, we make the following substitution of the $2n$ variables
$\{ \xvar{1}{i},\xvar{2}{i} \}_{i=1}^n$:
	\begin{equation} \label{substitute}
	\begin{aligned}
	\xvar{1}{i} &= x_i^{-\frac{1}{2}},\\
	\xvar{2}{i} &= x_i^\frac{1}{2}.
	\end{aligned}
	\end{equation}
The substitution \eqref{substitute} is harmless, because it is immediate
from
\eqref{decoupled-monomial} that the polynomial $f_{2,\ldots,2}(q,x)$ is
homogeneous of
total degree $2(2^n-1)$ in each of the sets of two variables $\{
\xvar{1}{i},\xvar{2}{i}
\}$.  Our result may now be stated as follows:
\begin{theorem}
\label{cube-theorem}
$$
\left[
   \sum_{T \in \Tree(Q_n)} q^{\dir(T)} x^{\dd(T)}
\right]_{\xvar{1}{i} \:=\: x_i^{-\frac{1}{2}}, \; \xvar{2}{i} \:=\:
x_i^\frac{1}{2}}
~=~ q_1 \cdots q_n \prod_{ \condns{A \subset [n]}{|A| \geq 2} }
          ~~ \sum_{i \in A} q_i \left( x_i^{-1} + x_i \right).
$$
\end{theorem}

\noindent
This result may shed light on the problem of finding a bijective proof
for the known number of spanning trees in the $n$-cube (see \cite[pp.
61-62]{Sta99b}).

Section~\ref{threshold-section} proves a result
(Theorem~\ref{threshold-factorization} below),
generalizing the Cayley-Pr\"ufer Theorem in two somewhat different
directions.  In one direction, it deals with threshold graphs, a
well-behaved
generalization of complete graphs.  Threshold graphs have many equivalent
definitions (see,
e.g., \cite[Chapters 7-8]{Merris-book}), but one that is convenient for
our purpose is the
following.  A graph $G$ is {\it threshold} if, after labelling its
vertices by
$[n]:=\{1,2,\ldots,n\}$ in weakly decreasing order of their degrees, the
degree sequence
$\lambda = (\lambda_1 \geq \cdots \geq \lambda_n)$ determines the graph
completely by the
rule that the neighbors of vertex $i$ are the $\lambda_i$ smallest members
of $[n]$ other
than $i$ itself. A result of Merris \cite{Merris} implies the following
generalization of
Cayley's formula to all threshold graphs.  It uses the notion of the {\it
conjugate
partition} $\lambda'$ to the degree sequence $\lambda$, whose Ferrers
diagram is obtained
from that of $\lambda$ by flipping across the diagonal.

\vskip .1in
\noindent
{\bf Merris' Theorem.}
{\it
Let $G$ be a threshold graph with vertices $[n]$ and degree sequence
$\lambda$.
Then the number of spanning trees in $G$ is $\prod_{r=2}^{n-1}
\lambda_r'$.
}
\vskip .1in

The natural vertex-ordering by degree for a threshold
graph $G$ induces a canonical edge orientation in any spanning tree $T$ of
$G$, by orienting
the edge $\{i,j\}$ from $j$ to $i$ if $j > i$.  Thus given a spanning tree
$T$ and a vertex
$i$,
one can speak of its {\it indegree} $\indeg_T(i)$ and {\it outdegree}
$\outdeg_T(i)$.

\begin{theorem}
\label{threshold-factorization}
Let $G$ be a connected threshold graph with vertices $[n]$ and degree sequence
$\lambda$.  Then
$$
\sum_{T \in \Tree(G)} \quad
   \prod_{i=1}^n x_i^{\indeg_T(i)}  y_i^{\outdeg_T(i)}
 ~=~  x_1 y_n \prod_{r=2}^{n-1}
    \left( \sum_{i=1}^{\lambda_r'} x_{\min\{i,r\}} y_{\max\{i,r\}}
\right).
$$

In particular, setting $y_i=x_i$ gives
$$
\sum_{T \in \Tree(G)} x^{\deg(T)} =
x_1 x_2 \cdots x_n \prod_{r=2}^{n-1} \left( \sum_{i=1}^{\lambda_r'} x_i
\right).
$$
\end{theorem}

\noindent
The proof, sketched in Section~\ref{threshold-section}, proceeds
by identification of factors.  The authors thank M. Rubey and an anonymous
referee for pointing out that it can also be deduced bijectively from a
very special case of a recent encoding theorem of Remmel and Williamson~\cite{RW}.

\section{Proof of Cayley-Pr\"ufer Theorem: the model}
\label{model-proof-section}

The goal of this section is to review Kirchhoff's Matrix-Tree Theorem, and use it to give a
proof of the Cayley-Pr\"ufer Theorem. Although this proof is surely known, we included it
both because we were unable to find it in the literature, and because it will serve as a
model for our other proofs.

Introduce a variable $e_{ij}$ for each edge $\{i,j\}$ in the complete graph $K_n$, with the
conventions that $e_{ij}=e_{ji}$ and $e_{ii}=0$. Let $L$ be the $n \times n$ weighted
Laplacian matrix defined by
	\begin{equation} \label{Laplacian-entries}
	L_{ij} := 
	\begin{cases}
	\sum_{k=1}^n e_{ik} & \text{ for } i=j \\
	- e_{ij} & \text{ for } i \neq j.
	\end{cases}
	\end{equation}

\vskip .1in
\noindent
{\bf Kirchhoff's Matrix-Tree Theorem}~\cite[\S 5.3]{Moon}
{\it
For any $r,s \in [n]$,
$$
\sum_{T \in \Tree(K_n)} ~~ \prod_{\{i,j\} \in T} e_{ij}
  = (-1)^{r+s} \det \hat{L}
$$
where $\hat{L}$ is the reduced Laplacian matrix obtained from $L$ by
removing row $r$ and column $s$.
}

\vskip .1in

We now restate and prove the Cayley-Pr\"ufer Theorem.

\vskip .1in
\noindent
{\bf Cayley-Pr\"ufer Theorem.}
{\it
$$
\sum_{T \in \Tree(K_n)} x^{\deg(T)}
~=~ x_1 x_2 \cdots x_n (x_1 + \cdots + x_n)^{n-2}
$$
where $x^{\deg(T)} := \prod_{i=1}^n x_i^{\deg_T(i)}$.
}
\begin{proof}
Apply the substitution $e_{ij}=x_i x_j$ to the weighted Laplacian matrix $L$
in Kirchhoff's Theorem.  Setting 
$f:=x_1 + \cdots + x_n$, one has from \eqref{Laplacian-entries}
$$
L_{ij} =
\begin{cases} 
(f - x_i) x_j & \text{ for } i=j \\
- x_i x_j     & \text{ for } i\neq j .
\end{cases}
$$
By Kirchhoff's Theorem, the left-hand side of the Cayley-Pr\"ufer Theorem coincides with the
determinant $\det\hat{L}$, where $\hat{L}$ is the reduced Laplacian $\hat{L}$ obtained from
this substituted $L$ by removing the last row and column. We wish to show that this
determinant coincides with the right-hand-side of the Cayley-Pr\"ufer Theorem.  Note that both
sides are polynomials in the $x_i$ of degree $2n-2$, and both have coefficient $1$ in
the monomial $x_1^{n-1} x_2 x_3 \cdots x_n$.  Therefore it suffices to show that the
determinant is divisible by
each of the variables $x_j$, and also by $f^{n-2}$. Divisibility by $x_j$ is clear since
$x_j$ divides every entry in the $j^{th}$ column of $L$ (and hence also $\hat{L}$).
Divisibility by $f$ follows from Lemma~\ref{identification-of-factors} below, once one
notices that in the quotient ring $\QQ[x_1,\ldots,x_n]/(f)$, this weighted Laplacian $L$
(and hence $\hat{L}$) reduces to a rank one matrix of the form $L = - v^T \cdot v$, where
$v$ is the row vector $\begin{bmatrix} x_1 & \cdots & x_n \end{bmatrix}$. \end{proof}

The following lemma, used in the preceding proof, is one of our main tools.  It generalizes
from one to several variables the usual statement on identification of factors in
determinants over polynomial rings (see \cite[\S 2.4]{Krattenthaler-adc}).
	
\begin{lemma}(Identification of factors)
\label{identification-of-factors}
Let $R$ be a Noetherian integral domain (e.g., a polynomial or Laurent ring in
finitely many variables over a field).  Let $f \in R$ be a prime element,
so that the quotient ring $R/(f)$ is an integral domain, and let
$K$ denote the field of fractions of $R/(f)$.
Let $A \in R^{n \times n}$ be a square matrix.  If the reduction 
$\bar{A} \in (R/(f))^{n \times n}$ has $K$-nullspace of dimension
at least $d$, then $f^d$ divides $\det A$ in $R$.
\end{lemma}

\begin{proof}
Let $\{\bar{v}_i\}_{i=1}^d$ be $d$ linearly independent vectors in $K^n$ lying in the
nullspace of $\bar{A}$.  Extend them to a basis $\{\bar{v}_i\}_{i=1}^n$ of $K^n$. By
clearing denominators, one may assume that $\{\bar{v}_i\}_{i=1}^n$ lie in $R/(f)^n$, and
then choose pre-images $\{v_i\}_{i=1}^n$ in $R$.

Letting $F$ denote the fraction field of $R$, we claim that $\{v_i\}_{i=1}^n$
is a basis for $F^n$.  To see this, assume not, so that there are
scalars $c_i \in F$ which are not all zero satisfying
	\begin{equation} \label{syzygy}
	\sum_{i=1}^n c_i v_i = 0.
	\end{equation}
Clearing denominators, one may assume that $c_i \in R$ for all $i$. If every $c_i$ is
divisible by $f$, one may divide the equation \eqref{syzygy} through by $f$, and repeat
this division until at least one of the $c_i$ is not divisible by $f$.  (This will happen
after finitely many divisions because $R$ is Noetherian.) But then reducing
\eqref{syzygy} modulo $(f)$ leads
to a nontrivial $K$-linear dependence among the vectors $\{\bar{v}_i\}_{i=1}^n$,
a contradiction.

Let $P \in F^{n \times n}$ be the matrix whose columns are the vectors $v_i$. Note that
$\det P$ is not divisible by $f$, or else the reductions $\{\bar{v}_i\}_{i=1}^n$ would
not form a $K$-basis in $K^n$.  Therefore, by Cramer's Rule, every entry of $P^{-1}$
belongs to the localization $R_{(f)}$ at the prime ideal $(f)$. Note the following
commutative diagram in which horizontal maps are inclusions and vertical maps are
reductions modulo $(f)$:
$$
\begin{CD}
R    &     @>>>  & R_{(f)} & @>>>  & F \\
@VVV &           & @VVV    &           \\
R/(f)&     @>>>  & K       & 
\end{CD}
$$

\noindent
Since $P^{-1}$ has entries in $R_{(f)}$, so does $P^{-1}AP$.
For each $i \in [d]$, the reduction of $Av_i$ vanishes in $K$, 
so every entry in the first $d$ columns
of $P^{-1}AP$ lies in the ideal $(f)$.  Hence $\det P^{-1}AP = \det A$ is divisible
by $f^d$ in $R_{(f)}$, thus also in $R$.
\end{proof}

The authors thank W.~Messing for pointing out a more general result, 
deducible by a variation of the above proof that uses Nakayama's Lemma:
\begin{lemma}
\label{Messing's-lemma}
Let $S$ be a (not necessarily Noetherian) local ring, with maximal ideal
$\mm$ and residue field $K:=S/\mm$.  Let $A \in S^{n \times n}$ be
a square matrix such that the reduction $\bar{A}$ has $K$-nullspace of dimension
at least $d$. Then $\det A \in \mm^d$.
\end{lemma}

Lemma~\ref{identification-of-factors} follows from this
by taking $S$ to be the localization $R_{(f)}$.

\section{Proof of Theorem~\ref{directions-theorem}}
\label{Cartesian-products-section}

We recall the statement of Theorem~\ref{directions-theorem}.

\vskip .1in
\noindent
{\bf Theorem~\ref{directions-theorem}.}
$$
\sum_{T \in \Tree(K_{n_1} \times \cdots \times K_{n_r})} 
                        q^{\dir(T)}
= \frac{1}{n_1 \cdots n_r} \prod_{\emptyset \neq A \subset [r]} 
    \left( \sum_{i \in A} q_i n_i \right )^{\prod_{i \in A}(n_i-1)}.
$$

As a prelude to the proof, we discuss some generalities about Laplacians and
eigenvalues of Cartesian products of graphs.  We should emphasize that
all results in this section refer only to {\it unweighted} Laplacians, that is,
one substitutes $e_{ij}=1$ for $i \neq j$ in the usual weighted
Laplacian $L(G)$ defined in \eqref{Laplacian-entries}.

The {\it Cartesian product} $G_1 \times \cdots \times G_r$ of graphs $G_i$ with
vertex sets $V(G_i)$ and edge sets $E(G_i)$ is defined as the graph with vertex set
$$
V(G_1 \times \cdots \times G_r) \;=\; V(G_1) \times \cdots \times V(G_r)
$$
and edge set
$$
\begin{aligned}
&E(G_1 \times \cdots \times G_r) \\
&\qquad \qquad =\;
\bigsqcup_{i=1}^r
V(G_1) \times \cdots \times V(G_{i-1}) \times E(G_i) \times  
  V(G_{i+1}) \times \cdots \times V(G_r).
\end{aligned}
$$
where $\bigsqcup$ denotes a disjoint union.  The following proposition follows easily from
this description; we omit the proof.

\begin{proposition}
\label{product-eigenvalues}
If $G_1,\ldots,G_r$ are graphs with (unweighted) Laplacian matrices $L(G_i)$, then
$$
L(G_1 \times \cdots \times G_r) =
\sum_{i=1}^r \id \otimes \cdots \id \otimes L(G_i) \otimes \id \otimes \cdots \id
$$
where $\id$ denotes the identity, and $L(G_i)$ appears in the $i^{th}$
tensor position.  

As a consequence, a complete set of eigenvectors for $L(G_1 \times \cdots \times G_r)$ can
be chosen of the form $v_1 \otimes \cdots \otimes v_r$, where $v_i$ is an eigenvector for
$L(G_i)$.  Furthermore, this eigenvector will have eigenvalue $\lambda_1 + \cdots +
\lambda_r$ if $v_i$ has eigenvalue $\lambda_i$ for $L(G_i)$.
\end{proposition}

We also will make use of the following variation of the Matrix-Tree
Theorem; see, e.g., \cite[Theorem~5.6.8]{Sta99b}.

\begin{theorem}
\label{Kirchhoff-variant}
If the (unweighted)
Laplacian matrix $L(G)$ has eigenvalues $\lambda_1,\ldots,\lambda_n$,
indexed so that $\lambda_n=0$, then the number of spanning trees in $G$ is
$$
\frac{1}{n} \lambda_1 \cdots \lambda_{n-1}.
$$
\end{theorem}

\noindent
{\it Proof of Theorem~\ref{directions-theorem}.} Both sides in the theorem are polynomials
in the $q_i$, and hence it suffices to show that they coincide whenever the $q_i$ are all
positive integers.  In that case, the left-hand side of the theorem has the following
interpretation.  Let $K_n^{(q)}$ denote the multigraph on vertex set $[n]$ having $q$
parallel copies of the edge $\{i,j\}$ for every pair of vertices $i,j$.  Then the left-hand
side of Theorem~\ref{directions-theorem} counts the number of spanning trees in the
Cartesian product
$$
K_{n_1}^{(q_1)} \x \cdots \x K_{n_r}^{(q_r)},
$$
as each spanning tree $T$ in $K_{n_1} \x \cdots \x K_{n_r}$ gives rise in an obvious way to
exactly $q^{\dir(T)}$ spanning trees in $K_{n_1}^{(q_1)} \x \cdots \x K_{n_r}^{(q_r)}$. It is
well-known that the (unweighted) Laplacian $L(K_n)$ has eigenvalues
$n, 0$ with multiplicities $n-1, 1$, respectively~\cite[Example~5.6.9]{Sta99b}.  Hence
$L(K_n^{(q)}) = q L(K_n)$ has eigenvalues $qn, 0$ with multiplicities $n-1, 1$,
respectively.  By
Proposition~\ref{product-eigenvalues}, $L(K_{n_1}^{(q_1)} \x \cdots \x K_{n_r}^{(q_r)})$ has
an eigenvalue $\sum_{i \in A} q_i n_i$ for each subset $A \subset [r]$, and this eigenvalue
occurs with multiplicity $\prod_{i \in A}(n_i-1)$.  As the zero eigenvalue arises (with
multiplicity $1$) only by taking $A = \emptyset$, Theorem~\ref{Kirchhoff-variant} implies
that the number of spanning trees in $K_{n_1}^{(q_1)} \x \cdots \x K_{n_r}^{(q_r)}$ is
$$
\frac{1}{n_1 \cdots n_r} \prod_{\emptyset \neq A \subset [r]} 
   \left( \sum_{i \in A} q_i n_i \right)^{\prod_{i \in A}(n_i-1)}.
\eqno\qed
$$

\section{Proof of Theorem~\ref{divisibility-theorem}}
\label{divisibility-section}

We recall here the statement of Theorem~\ref{divisibility-theorem}.

\noindent
{\bf Theorem~\ref{divisibility-theorem}.}
{\it The spanning tree enumerator
$$
f_{n_1,\ldots,n_r}(q,x):=
\sum_{T \in \Tree(K_{n_1} \times \cdots \times K_{n_r})} 
                        q^{\dir(T)} x^{\dd(T)}
$$
is divisible by 
	$$q_i^{n_i-1},$$
by
	$$\left( \xvar{j}{i} \right)^{n_1 \cdots n_{i-1} n_{i+1} \cdots n_r}$$
and by
	$$\left( \xvar{1}{i} + \cdots + \xvar{n_i}{i} \right)^{n_i-2}$$
for each $i \in [r]$ and $j \in [n_i]$.}

\begin{proof}
To see divisibility by $q_i^{n_i-1}$, note that every spanning tree in $K_{n_1} \times \cdots
\times K_{n_r}$ is connected, and hence gives rise to a connected subgraph of $K_{n_i}$ when
one contracts out all edges not lying in direction $i$.  This requires at least $n_i-1$
edges in direction $i$ in the original tree.

To see divisibility by $(\xvar{j}{i})^{n_1 \cdots n_{i-1} n_{i+1} \cdots n_r}$, note that every
spanning tree has an edge incident to each vertex, and therefore to each of the $n_1 \cdots
n_{i-1} n_{i+1} \cdots n_r$ different vertices which have $i^{th}$ coordinate equal to some
fixed value $j \in [n_i]$.

Lastly, we check divisibility by $(\xvar{1}{i} + \cdots + \xvar{n_i}{i})^{n_i-2}$.
Starting with the weighted Laplacian matrix~(\ref{Laplacian-entries})
for $K_{n_1} \times \cdots \times K_{n_r}$ (regarded as a subgraph of
$K_{n_1\cdots n_r}$), let $L$ be the matrix obtained by
the following substitution: if
$\{k,l\}$ represents an edge of $K_{n_1} \times \cdots \times K_{n_r}$ in direction $i$
between the two vertices $k=(k_1,\ldots,k_r)$ and $l=(l_1,\ldots,l_r)$, then we set
	\begin{equation} \label{ekl-substitution}
	e_{kl}=
	q_i \cdot
	\xvar{k_1}{1} \cdots \xvar{k_r}{r} \cdot
	\xvar{l_1}{1} \cdots \xvar{l_r}{r}\;;
	\end{equation}
otherwise we set $e_{kl}=0$.
Then Kirchhoff's Theorem says that $f_{n_1,\ldots,n_r}(q,x) = \pm \det{\hat{L}}$
for any reduced matrix $\hat{L}$ obtained from $L$ by removing a
row and column.  Thus by Lemma~\ref{identification-of-factors}, it suffices
for us to show that $\hat{L}$ has nullspace of dimension at least
$n_i-2$ modulo $f^{(i)} := \xvar{1}{i} + \cdots + \xvar{n_i}{i}$.
In fact, we will show that $L$ itself has nullspace of dimension at
least $n_i-1$ in this quotient.  To see this, one can check that, as
in Proposition~\ref{product-eigenvalues}, the matrix $L$ has the following
simpler description, due to our ``decoupling'' substitution of variables:
$$
L = \sum_{i=1}^r X^{(1)} \otimes \cdots \otimes X^{(i-1)} 
       \otimes q_i L^{(i)} \otimes X^{(i+1)} \otimes \cdots \otimes X^{(r)},
$$
where $X^{(i)}$ is the diagonal matrix with entries 
$(\xvar{1}{i})^2,\ldots,(\xvar{n_i}{i})^2$, and $L^{(i)}$ is
obtained by making the substitution $e_{kl} = \xvar{k}{i} \xvar{l}{i}$
in the weighted Laplacian matrix for $K_{n_i}$.
In the proof of the Cayley-Pr\"ufer Theorem,
we saw that $L^{(i)}$ has rank~$1$ modulo $(f^{(i)})$, and thus a nullspace
of dimension $n_i-1$.  If $v$ is any nullvector for $L^{(i)}$ modulo $(f^{(i)})$,
then the following vector is a nullvector for $L$ modulo $(f^{(i)})$:
$$
{\bf 1}_{n_1} \otimes \cdots \otimes {\bf 1}_{n_{i-1}} 
       \otimes v \otimes {\bf 1}_{n_{i+1}} \otimes \cdots \otimes {\bf 1}_{n_r}
$$
where ${\bf 1}_{m}$ represents a vector of length $m$ with all entries equal to $1$.
Since varying $v$ leads to $n_i-1$ linearly independent
such nullvectors, the proof is complete.
\end{proof}

\section{Proof of Theorem~\ref{cube-theorem}}
\label{cube-section}

We recall here the statement of Theorem~\ref{cube-theorem}, using slightly different
notation.  Regard the vertex set of $Q_n$ as the power set $2^{[n]}$, so that vertices
correspond to subsets of $S \subset [n]$. For any subset $S \subset [n]$, let $x_S :=
\prod_{i \in S} x_i$. Write $x^{\ddq(T)}$ for the decoupled degree monomial corresponding to
a tree $T$ under the substitution \eqref{substitute}, that is,
	\begin{equation} \label{decoupled-monomial-for-cube}
	\begin{aligned}
	x^{\ddq(T)} ~&=~ \prod_{S \subset [n]}
		\left( \frac{x_S}{x_{[n] \sm S}} \right)^{\frac{1}{2} \deg_T(S)} \\
	~&=~ \prod_{\text{ edges }\{S,R\} \text{ in }T} \frac{x_S x_R}{x_{[n]}}
	\end{aligned}
	\end{equation}

\noindent
{\bf Theorem~\ref{cube-theorem}.}
$$
\sum_{T \in \Tree(Q_n)} q^{\dir(T)} x^{\ddq(T)} 
~=~ q_1 \cdots q_n \prod_{\condns{A \subset [n]}{|A| \geq 2}} \quad
	\sum_{i \in A} q_i \left( x_i^{-1} + x_i \right).
$$

\begin{proof}
As before, regard the vertex set of $Q_n$ as the power set $2^{[n]}$.  Denote the symmetric
difference of two sets $S$ and $R$ by $S \sd R$, and abbreviate $S \sd \{i\}$ by $S \sd i$.
Thus two vertices $S,R$ form an edge in $Q_n$ exactly when $S \sd R$ is a singleton set
$\{i\}$; in this case the \defterm{direction} of this edge is $\dir(e):=i$.  It is useful to
note that the neighbors of $S$ are
	$$N(S) = \{ S \sd i \st i \in [n] \}.$$
Our goal is to show that the two sides of the theorem coincide as elements of $\ZZ[q_i, x_i,
x_{i}^{-1}]$. Note that the two sides coincide as polynomials in $q_i$ after setting $x_i=1$
for all $i$, using the special case of Theorem~\ref{directions-theorem} in which all
$n_i=2$.

We next show that both sides have the same maximum and minimum total degrees as Laurent
polynomials in the $x_i$.  Each side is easily seen to be invariant under the substitution
$x_i \mapsto x_i^{-1}$ (this follows from the antipodal symmetry of the $n$-cube for the
left-hand side), so it suffices to show that both sides have the same maximum total degree.
For the right-hand side, the maximum total degree in the $x_i$ is simply the number of
subsets $S \subset [n]$ with $|S| \geq 2$, that is, $2^n-n-1$.

For the left-hand side, we argue as follows.  Denote by $V'$ the set of vertices of $Q_n$
other than $[n]$.  For any spanning tree $T$ and vertex $S \in V'$, define $\phi(S)$ to be
the parent vertex of the vertex $S$ when the tree $T$ is rooted at the vertex $[n]$ (that
is, $\phi(S)$ is the first vertex on the unique path in $T$ from $S$ to $[n]$);  then the
edges of $T$ are precisely
	$$E(T) = \left\{ \{S,\phi(S)\} \st S \in V' \right\}.$$
One has $|\phi(S)| = |S| \pm 1$ (because $\phi(S) = S \sd i$ for some $i$).  Therefore the
total degree of the monomial $x^{\ddq(T)}$ will be maximized when $|\phi(S)|=|S|+1$ for all $S
\in V'$; for instance, when $\phi(S) = S \cup \{\max([n] \sm S)\}$.  In this case, that
total degree is
	$$
	\sum_{\{S,\phi(S)\} \in T} \left( |S|+|\phi(S)|-n \right) \ = \
	\sum_{S \in V'} \left( 2|S|+1-n \right) \ = \ 
	2^n-n-1.
	$$

Having shown that both sides have the same total degree in the $q_i$, the same maximum and
minimum total degrees in the $x_i$, and that they coincide when all $x_i=1$, it suffices by
unique factorization to show that the left-hand side is divisible by each factor on the
right-hand side, that is, by
$$
f_A ~:=~ \sum_{i \in A} q_i \left( x_i^{-1}+x_i \right).
$$
Henceforth, fix $A \subset [n]$ of cardinality $\geq 2$.  It is not hard to check that 
$f_A$ is irreducible in $\ZZ[q_i, x_i, x^{-1}_i]$, using the fact that it is a linear form 
in the $q_i$.

Starting with the weighted Laplacian matrix~(\ref{Laplacian-entries}) for $Q_n$,
whose rows and columns are indexed by subsets $S \subseteq [n]$,
let $L$ be the matrix obtained by making the substitutions
	$$e_{S,R} = \begin{cases}
	\displaystyle \frac{q_i x_S x_{S \sd i}}{x_{[n]}} & \text{for } S \sd R = \{i\}, \\
	0 & \text{for } |S \sd R| > 1.
	\end{cases}$$
By Kirchhoff's Theorem, the left-hand side in Theorem~\ref{cube-theorem} is the determinant
of the reduced Laplacian matrix $\hat{L}$ obtained from $L$ by removing the row and column
indexed by $S=\emptyset$.  It therefore suffices to show that the reduction of $\hat{L}$
modulo $(f_A)$ has nontrivial nullspace.  We will show that
	\begin{equation} \label{nullvector}
	v := \sum_{\emptyset \neq S \subset [n]}
	\left( x_A^2 - (-1)^{|A \cap S|} (x_{A \sm S})^2 \right) \ebv_S
	\end{equation}
is a nullvector\footnote{The form of the nullvector~(\ref{nullvector}) was suggested by 
experimentation using the computational commutative algebra package {\tt 
Macaulay}~\cite{Macaulay} to compute the nullspace of $\hat{L}$ in the quotient ring modulo 
$(f_A)$.}, where $\ebv_S$ is the standard basis vector corresponding to $S$.
Note that the entries of $v$ are not all zero modulo $(f_A)$; it remains to
check that every entry $(\hat{L}v)_R$ of $\hat{L}v$ is a multiple of $f_A$.
Since $\hat{L}_{R,S}=0$ unless $S=R$ or $S=R \sd i$ for some $i$, one has
	\begin{align}
	(\hat{L}v)_R ~&=~ \hat{L}_{R,R} \: v_R ~+~ \sum_{i=1}^n \hat{L}_{R,R \sd i} \:
		v_{R \sd i} \notag \\
	&=~ \sum_{i=1}^n \frac{q_i x_R x_{R \sd i}}{x_{[n]}} \left( x_A^2 - (-1)^{|A \cap R|}
		(x_{A \sm R})^2 \right) \notag \\
	&\qquad \qquad -~ \sum_{i=1}^n \frac{q_i x_R x_{R \sd i}}{x_{[n]}} \left( x_A^2 -
		(-1)^{|A \cap (R \sd i)|} (x_{A \sm (R \sd i)})^2 \right) \notag \\
	\label{bigmess} &=~ \frac{x_R}{x_{[n]}} \sum_{i=1}^n q_i x_{R \sd i} \left(
	(-1)^{|A \cap (R \sd i)|} (x_{A \sm (R \sd i)})^2 - (-1)^{|A \cap R|}
	(x_{A \sm R})^2 \right)~.
	\end{align}
If $i \not\in A$, then $A \cap R = A \cap (R \sd i)$ and $A \sm R = A \sm (R \sd i)$, so the
summand in \eqref{bigmess} is zero.  If $i \in A$, then $|A \cap R| = |A \cap (R \sd i)| \pm
1$, so one may rewrite \eqref{bigmess} as follows:
	\begin{equation} \label{bigmess2}
	(\hat{L}v)_R ~=~ -(-1)^{|A \cap R|} \; \frac{x_R}{x_{[n]}} \sum_{i \in A}
	q_i x_{R \sd i} \left(	(x_{A \sm (R \sd i)})^2 + (x_{A \sm R})^2 \right).
	\end{equation}
Note also that when $i \in A$,
	\begin{equation*}
	x_{R \sd i} \;=\; \begin{cases}
		x_R \, x_i^{-1} &\text{for } i \in R \\
		x_R \, x_i &\text{for } i \not\in R
	\end{cases}
	\end{equation*}
and
        \begin{equation*}
	x_{A \sm (R \sd i)} \;=\; \begin{cases}
		x_{A \sm R} \, x_i &\text{for } i \in R \\
		x_{A \sm R} \, x_i^{-1} &\text{for } i \not\in R.
	\end{cases}
	\end{equation*}
Therefore one may rewrite \eqref{bigmess2} as follows:
	\begin{eqnarray*}
	(\hat{L}v)_R  &=&
	    \pm\frac{x_R}{x_{[n]}} \left( \sum_{i \in A \cap R} q_i x_R x_i^{-1} \left(
	    (x_{A \sm R} \, x_i)^2 + (x_{A \sm R})^2 \right) \right. \\
	& & \left. \qquad +~ \sum_{i \in A \sm R} q_i x_R x_i \left( (x_{A \sm R} \, 
	    x_i^{-1})^2 + (x_{A \sm R})^2 \right) \right) \\
	&=& \pm\frac{(x_R x_{A \sm R})^2}{x_{[n]}} \left( \sum_{i \in A \cap R}
	    q_i x_i^{-1} (x_i^2+1) ~+~ \sum_{i \in A \sm R} q_i x_i (x_i^{-2}+1) \right) \\
	&=& \pm\frac{(x_R x_{A \sm R})^2}{x_{[n]}} \sum_{i \in A} q_i(x_i+x_i^{-1}) \\
	&=& \pm\frac{(x_R x_{A \sm R})^2}{x_{[n]}} f_A,
	\end{eqnarray*}
which shows that $(\hat{L}v)_R$ is zero modulo $(f_A)$ as desired.
\end{proof}

\section{Proof of Theorem~\ref{threshold-factorization}}
\label{threshold-section}

We recall the statement of Theorem~\ref{threshold-factorization}.

\vskip .1in
\noindent
{\bf Theorem~\ref{threshold-factorization}.}
{\it 
Let $G$ be a connected threshold graph with vertices $[n]$, edges $E$, and degree sequence 
$\lambda$.  Then 
$$
\sum_{T \in \Tree(G)} \quad
   \prod_{i=1}^n x_i^{\indeg_T(i)}  y_i^{\outdeg_T(i)} 
 ~=~  x_1 y_n \prod_{r=2}^{n-1} 
    \left( \sum_{i=1}^{\lambda_r'} x_{\min\{i,r\}} y_{\max\{i,r\}} \right).
$$
}

As noted in the Introduction, this result is a special case of Theorem~2.4 of~\cite{RW}.
For this reason, and because the ideas of the proof are quite similar to those of 
Theorem~\ref{cube-theorem}, we omit most of the technical details.  We write $N(v)$ for the 
neighbors of a vertex $v$, and denote the set $\{i,i+1,\ldots,j\}$ by $[i,j]$.

{\it Sketch of proof.\/}
The partitions $\lambda$ which arise as degree sequences of threshold graphs have been
completely characterized (see, e.g.,~\cite[Theorem 8.5]{Merris-book}).  In particular, suppose 
that
Durfee square of $\lambda$ (the largest square which is a subshape of $\lambda$) has side 
length~$s$.  Then for all $r \in [n]$,
	\begin{equation} \label{useful-eqn}
	\begin{array}{rl}
	\text{either} \ & r \leq s < \lambda'_r \;=\; 1+\lambda_r \\
	\text{or}     \ & r > s \geq \lambda'_r \;=\; \lambda_{r+1}.
	\end{array}
	\end{equation}
Using these identities, one may rewrite the desired equality as
	\begin{equation} \label{threshold-rewrite}
	\sum_{T \in \Tree(G)} \ \ 
	   \prod_{i=1}^n ~ x_i^{\indeg_T(i)}  y_i^{\outdeg_T(i)} \\
	 ~=~ x_1 \cdot \prod_{r=2}^s f_r \cdot \prod_{r=s+1}^{n-1} g_r \cdot
		\prod_{r=s+1}^n y_r.
	\end{equation}
where
	$$
	f_r \ := \ y_r \sum_{i=1}^r x_i \ + \ x_r \sum_{i=r+1}^{1+\lambda_r} y_i
		\qquad \text{and} \qquad
	g_r \ := \ \sum_{i=1}^{\lambda_{r+1}} x_i.
	$$

Both the left-hand and right-hand sides of~(\ref{threshold-rewrite})
are polynomials in the $x_i, y_i$ of total degree
$2n-2$, and both have coefficient of $x_1^{n-1}y_2 y_3 \cdots y_n$ equal to $1$ (because $N(1) 
= [2,n]$).  Thus it
suffices to prove that the left-hand side is divisible by each of the factors on the
right-hand side. By Kirchhoff's Theorem, this left-hand side is the determinant of the matrix
$\hat{L}$ obtained from the usual weighted Laplacian matrix by removing the first row and
column and making the substitution
$$
e_{ij} = 
\begin{cases}
x_{\min\{i,j\}} y_{\max\{i,j\}}	& \text{ for } \{i,j\} \in E \\
0				& \text{ for } \{i,j\} \not\in E.
\end{cases}
$$

\underline{First}, one must show that the left-hand side of the theorem is divisible by 
the monomial factor $x_1 \prod_{r=s+1}^n y_r$.  Every spanning tree
$T$ of $G$ contains an edge of the form $\{1,j\}$, which contributes a factor of $x_1y_j$ to
the monomial corresponding to $T$.  In particular, $x_1$ divides the left-hand side
of~(\ref{threshold-rewrite}).  
Furthermore, if $r>s$, then $q<r$ whenever $q \in N(r)$.  In
particular, $y_r$ divides every 
entry in the $r^{th}$ row of $\hat{L}$.

\underline{Second}, one must show that $f_r$ divides $\det \hat{L}$ for $r \in [2,s]$. Clearly
$f_r$ is irreducible, since neither sum in the definition of $f_r$ is empty.  Define a column
vector\footnote{Found using {\tt Macaulay}; see the earlier footnote.}
	$$v = \sum_{i=1}^r x_i \ebv_r + \sum_{i=r+1}^{\lambda'_r} x_r \ebv_i,$$
where $\ebv_i$ denotes the $i^{th}$ standard basis vector. Note that the entries of $v$ are
not all divisible by $f_r$, so that $v$ is a non-zero vector modulo $(f_r)$.  By 
Lemma~\ref{identification-of-factors}, it is now sufficient to show that for each $j$,
the entry $(\hat{L}v)_j$ of $\hat{L}v$ is divisible by $f_r$.  One must consider four 
cases depending on the value of $j$: (i) $j < r$, (ii) $j=r$, (iii) $j > r$ and 
$\{j,r\} \in E$, (iv) $j > r$ and $\{j,r\} \notin E$.  We omit the
routine calculations, which are similar to the proof of Theorem~\ref{cube-theorem}.

\underline{Third}, one must show that $g_r$ divides $\det\hat{L}$ for all $r \in [s+1,n-1]$.
In fact, some higher power of $g_r$ may divide $\det\hat{L}$, as we now explain.
If $\lambda$ has exactly $b$ columns of height $\lambda'_r$, i.e.,
        $$\lambda'_{a-1} > \lambda'_a = \;\cdots\; = \lambda'_r = \;\cdots\; = \lambda'_{a+b-1}
        > \lambda'_{a+b}$$
for some $a>s$, then $N(i) = [\lambda_r]$ for all vertices $i \in [a+1,a+b]$, so
        $$g_a = g_{a+1}= \;\cdots\; = g_r = \;\cdots\; = g_{a+b-1}.$$
Accordingly, one must show that $g_r^b$ divides $\det\hat{L}$.  Restricting the Laplacian 
matrix $L$ to the columns $[a+1,a+b]$ yields a rank-$1$ matrix of the form
	$$
	-\begin{bmatrix} x_1 & x_2 & \cdots & x_{\lambda_r} & 0 & \cdots & 0 \end{bmatrix}^T 
	\cdot \begin{bmatrix} y_{a+1} & y_{a+2} & \cdots & y_{a+b} \end{bmatrix}.
	$$
Consequently, both $L$ and $\hat{L}$ have 
$b-1$ linearly independent nullvectors modulo $(g_r)$ supported in
coordinates $[a+1,a+b]$.  It remains only to exhibit one further nullvector 
for $\hat{L}$ which is supported in at least one
coordinate outside that range.
We claim that such a
vector\footnote{Found using {\tt Macaulay}; see the earlier footnote.} is
$$
   	v = \sum_{i=1+\lambda_r}^a \left( y_{a+b} \ebv_i - y_i \ebv_{a+b} \right).
$$
One must verify that for each $k$, the
$k^{th}$ coordinate $(\hat{L}v)_k$ vanishes modulo $(g_r)$.
This calculation splits into four cases:
(i) $k \in [2,\lambda_r]$,
(ii) $k \in [1+\lambda_r,a]$,
(iii) $k \in [a+1,n-1] \sm \{a+b\}$,
(iv) $k=a+b$.
Once again, we omit the routine verification.
\qed


\begin{thebibliography}{99}

\bibitem{Macaulay}
D. Bayer and M. Stillman,
Macaulay: A computer algebra system for algebraic geometry.
Software, 1994.
Available at {\tt www.math.columbia.edu/$\sim$bayer/Macaulay/}.

\bibitem{Krattenthaler-adc}
C. Krattenthaler,
Advanced determinant calculus.
The Andrews Festschrift (Maratea, 1998). 
{\it S\'em. Lothar. Combin.} {\bf 42} (1999), Art. B42q, 67 pp. (electronic).

\bibitem{Merris}
R. Merris, 
Degree maximal graphs are Laplacian integral.
{\it Linear Algebra Appl.} {\bf 199} (1994), 381--389.

\bibitem{Merris-book}
R. Merris, 
Graph theory. 
Wiley-Interscience Series in Discrete Mathematics and Optimization. 
Wiley-Interscience, New York, 2001.

\bibitem{Moon}
J.W. Moon, 
Counting labelled trees. 
{\it Canadian Mathematical Monographs} {\bf 1}
Canadian Mathematical Congress, Montreal, 1970.

\bibitem{RW}
J. Remmel and S.G. Williamson,
Spanning trees and function classes.
{\it Electron. J. Combin.} {\bf 9} (2002), Research Paper 34, 24 pp. (electronic).

\bibitem{Sta99b}
R.P. Stanley,
Enumerative Combinatorics, Volume II.
\textsl{Cambridge Studies in Advanced Mathematics} {\bf 62}.
Cambridge University, 1999.

\end{thebibliography}
\end{document}